\documentclass[12pt]{article}
\usepackage{amsmath,amssymb,amsfonts,epsfig,xspace,picinpar}
\newcommand{\Cov}{\operatorname{Cov}}
\newcommand{\Var}{\operatorname{Var}}

\begin{document}
\pagestyle{plain}
\thispagestyle{empty}

\title{\vspace*{-50pt}
Diagonal Sums of Boxed Plane Partitions}
\author{
\begin{tabular}{c}
David B. Wilson \\[-2pt]
 \small Microsoft Research 
\end{tabular}
}
\date{}

\maketitle

\vspace*{-20pt}
\begin{quote}\textbf{Abstract:}
We give a simple proof of a nice formula for the means and covariances of the
diagonal sums of a uniformly random boxed plane parition.
\end{quote}

\begin{figwindow}[0,r,%
{\label{fig:cube}\epsfig{figure=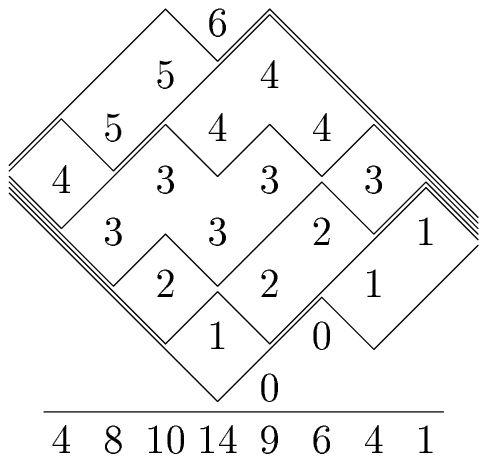}},%
{A $4\times5\times 6$ boxed plane partition with its contours and diagonal sums.}]
An $a\times b\times c$ boxed plane partition is an $a\times b$ grid of
integers between $0$ and $c$ inclusive, such that the numbers decrease
weakly in each row and column.  At the right is a $4\times 5\times
6$ boxed plane parition, which for convenience we have drawn rotated $45^\circ$.
We have added up these numbers in the direction along the main
diagonal of the $a\times b$ box to obtain the diagonal sums
$S_{-a+1},\ldots,S_{b-1}$.  If we pick the boxed plane partition
uniformly at random, these form a sequence of random variables, and we
show that their means and covariances are given by
\end{figwindow}
$$ E[S_i] = \begin{cases}(a+i) bc/(a+b) & i\leq 0\\ (b-i) ac/(a+b) & i\geq 0\end{cases} \hspace{2cm} $$
$$ \Cov(S_i,S_j) = (a+i) (b-j) \times \frac{abc(a+b+c)}{(a+b)^2((a+b)^2-1)}\ \ \ \ (i\leq j).$$
Notice that while the expected values ``see'' the corner at the origin,
the other corner does not enter the formula, and neither corner enters the formula for the covariances.
For given values of $a$, $b$, and $c$, the covariances are just proportional
to the product of the distances from $i$ and $j$ to the endpoints.
As Kenyon points out, a similar covariance property holds for Brownian
bridges, and indeed can be deduced from this formula by taking $c=1$
and $a,b\rightarrow\infty$.  We are unaware of similarly nice formulas
for e.g.\ the row sums.

As Stembridge points out,
diagonal sums appear in some generating functions such as
$$ \sum_{\text{$a\times b\times\infty$ bpp's}} \prod_i x_i^{S_i} = \prod_{i=-a+1}^0\prod_{j=0}^{b-1} \frac{1}{1-x_i x_{i+1}\cdots x_j},$$
which is due to Stanley (see \cite[Chapter 7]{EC2}).
The corresponding generating function for $a\times b\times c$ bpp's is
not so nice, but Krattenthaler \cite[pp 192]{K} expresses it in terms of a determinant.
We do not know a derivation of the covariance formula using this
approach.

To prove these formulas we look at the contours associated with a
boxed plane partition.  The contours come from viewing the numbers in
the plane partition as heights, so that between adjacent cells with
heights $z_1$ and $z_2$ there will be $|z_1-z_2|$ contours (see
figure).  The contours are noncrossing, but may share vertices and
edges.

Consider the $i$th diagonal, and condition on the locations where the
contours cross the two adjacent diagonals.  Each contour will either
make two down moves or two up moves, or it will be flexible and make an up and down
move in some random order.  The flexible lattice paths do not interact
at all unless they intersect diagonals $i-1$ and $i+1$ in the same
locations.  Since the boxed plane partition is uniformly random, if a
group of $k$ flexible contours start and end at the same locations,
the expected number that go down and then up is $k/2$.  If we
let $Y_i$ be the sum over contours of the height of the contour above
the line connecting the left and right corners
($Y_i+S_i$ is a deterministic function of $i$), then
we have $$E[Y_i|\text{diagonals $i-1$ and $i+1$}]=(Y_{i-1}+Y_{i+1})/2.$$

Since $Y_{-a}=0=Y_b$, we find $E[Y_i]=0$, which leads to the formula for $E[S_i]$.
Next observe
\begin{align*}
E[Y_i Y_j] &= \sum_y \Pr[Y_i=y] y E[Y_j|Y_i=y]
 = \sum_y \Pr[Y_i=y] y \frac{b-j}{b-i} y = \frac{b-j}{b-i} E[Y_i^2].
\intertext{Similarly}
E[Y_i Y_j] &= \sum_y \Pr[Y_j=y] y E[Y_i|Y_j=y]
 = \sum_y \Pr[Y_j=y] y \frac{a+i}{a+j} y = \frac{a+i}{a+j} E[Y_j^2].
\end{align*}
Equating these formulas gives
$E[Y_j^2]=(a+j)(b-j)/(a+b-1)E[Y_{a-1}^2]$, and hence
\begin{equation}
\Cov(S_i,S_j) = \Cov(Y_i,Y_j) = E[Y_i Y_j] = \frac{(a+i)(b-j)}{(a+b-1)} \Var[S_{-a+1}].\label{cvij}
\end{equation}
To compute $\Var[S_{-a+1}]$ we let $S=\sum_i S_i$ and write
\begin{align}
\Var[S] &= 2\sum_{i<j} \Cov(S_i,S_j) + \sum_{i} \Cov(S_i,S_i) \notag \\
        &= \left[2\sum_{-a<i<j<b} (a+i)(b-j) + \sum_{-a<i=j<b} (a+i)(b-j)\right]
           \frac{\Var[S_{-a+1}]}{a+b-1} \notag \\
        &= \frac{(a+b)^2((a+b)^2-1)}{12} \frac{\Var[S_{-a+1}]}{a+b-1} \label{vx1}
\end{align}
Next $\Var[S]$ can be computed using the $q$-analogue of MacMahon's formula
$$\sum_{\text{$a\times b\times c$ bpp's}} q^S = \prod_{i=1}^a \prod_{j=1}^b \prod_{k=1}^c \frac{(i+j+k-1)_q}{(i+j+k-2)_q},$$
where $n_q=q^{n-1}+\cdots+q+1$.  This calculation, perhaps first carried out by Blum \cite{B}, yields
\begin{equation}\Var[S]=abc(a+b+c)/12\label{vx2}.\end{equation}
Combining \eqref{cvij}, \eqref{vx1}, and \eqref{vx2} yields the formula for the covariance.

\bibliographystyle{plain}

\end{document}